\documentclass[preprint,review, 12pt,3p,times]{elsarticle}

\usepackage{amsthm,amsmath,amsfonts,graphicx}

\usepackage[english]{babel}
\usepackage{amsmath,amssymb}
\usepackage{amsfonts}
\usepackage{latexsym}
\usepackage{epsfig}
\usepackage{graphicx}
\usepackage{url}
\usepackage[T1]{fontenc}
\usepackage{color}

\newcommand{\ccr}{\color{red}}
\newcommand{\cb}{\color{blue}}
\definecolor{Dark}{gray}{.20}
\definecolor{Brown}{cmyk}{0, 0.8, 1, 0.6}
\definecolor{Yellow}{rgb}{1, 1, 0}
\definecolor{Light}{gray}{.80}
\definecolor{Black}{rgb}{0, 0, 0}
\newcommand{\cbr}{\color{Brown}}

\newtheorem{thm}{\ccr Theorem}
\newtheorem{lem}{\cbr Lemma}
\newtheorem{prop}{\cb  Proposition}
\newtheorem{cor}{\cb  Corollary}
\newtheorem{defi}{ Definition}
\newdefinition{exa}{Example}
\newdefinition{rem}{Remark}
\newdefinition{conj}{Conjecture}

\begin{document}

\begin{frontmatter}

\title{On the geodetic and the hull numbers in strong product graphs\tnoteref{t1}}


\author[ual]{J.~Cáceres}
\ead{jcaceres@ual.es}

\author[upc]{C.~Hernando}
\ead{carmen.hernando@upc.edu}

\author[upc]{M.~Mora}
\ead{merce.mora@upc.edu}

\author[upc]{I. M.~Pelayo}
\ead{ignacio.m.pelayo@upc.edu}

\author[ual]{M. L. Puertas}
\ead{mpuertas@ual.es}


\address[ual]{Universidad de Almería, Almería, Spain }

\address[upc]{Universitat Politècnica de Catalunya, Barcelona, Spain}

\author{}

\address{}

\begin{abstract}

A set $S$ of vertices  of a connected graph $G$ is convex, if for any pair of vertices  $u,v\in S$, every shortest path joining $u$ and $v$ is contained in $S$. The convex hull $CH(S)$ of a set of vertices $S$ is defined as the smallest convex set in G containing $S$. The set $S$ is  geodetic, if every vertex of $G$  lies on some shortest path joining two vertices in S, and it is said to be a hull set if  its convex hull is $V(G)$. The geodetic and the hull numbers of $G$ are the cardinality of a minimum geodetic and a minimum hull set, respectively.
In this work, we investigate the behavior of both  geodetic and  hull sets with respect to  the  strong product operation for graphs. We also stablish some bounds for  the geodetic number and the hull number and obtain the exact value of these  parameters for a number of strong product graphs.

\end{abstract}


\end{frontmatter}

\section{Introduction}
\label{int}

The process of rebuilding a network modelled by a connected graph is a discrete optimization problem, consisting in finding a subset of vertices of cardinalilty as small as possible, which, roughly speaking, would allow us to store and retrieve the whole graph.  One way to approach this problem is by using  a certain convex operator. This procedure  
has attracted much attention since
it was shown by Farber an Jamison \cite{FJ:86} that every convex
subset in a graph is the convex hull of its extreme vertices if an
only if the graph is chordal and contains no induced 3-fan. From then on, a number of variants of this  approach have been proposed   \cite{CMOP:05,P:08}. One of them, consists in using, instead of the convex hull operator, the closed interval operator, i.e., considering geodetic sets instead of hull sets \cite{CHMPPS:06,CHMPPS:08}. 
Unfortunately,  computing  geodetic sets and  hull sets of minimum cardinality, are known to be NP-hard problems for general graphs  \cite{DGKPS:09,DPRS:10}. This fact has motivated the study of these two problems for graph classes which can be obtained by means of graph operations,   such as cartesian product \cite{BKT:08,CC:04,JPP}, composition
\cite{CC:05} and join \cite{CCG:06}. Let us notice that in these graphs, infomation about factor graphs can be used to obtain geodetic and hull sets and to compute geodetic and hull numbers.

In this work, we
study both geodetic and hull sets of minimum cardinality, in strong product graphs. This graph operation has been extensively investigated in relation to a wide range of subjects, including: connectivity \cite{BS:07}, pancyclicity \cite{KMPS:04,KS:08},  chromaticity \cite{Z:06}, bandwidth \cite{K:08}, independency \cite{H:73,V:98}  and  primitivity \cite{K:93}.  Section 2 is devoted to introduce the main definitions and notation used throughout the paper. In Section 3, we study the behavior of  geodetic and hull sets with respect to the strong product operation. In Section 4, a number of lower and upper sharp bounds for the geodetic number and the hull number of the strong product of two graphs are presented. Finally, the last  ection is devoted to obtain the exact value of  the geodetic number and the hull number of the strong product of some basic families of graphs, such as paths, complete graphs and cycles.

\section{Graph theoretical preliminaries}

We consider only finite, simple, connected graphs. For undefined
basic concepts we refer the reader to introductory graph theoretical
literature, e.g., \cite{west}. Given vertices $u,v$ in a graph $G$
we let $d_G(u,v)$ denote the distance between $u$ and $v$ in $G$.
When there is no confusion, subscripts will be omitted. 
The diameter $diam(G)$ of $G$ is the maximum distance between any two vertices of G.
An $x-y$
path of length $d(x,y)$ is called an $x-y$ \emph{geodesic}.  The
\emph{closed interval} $I[x,y]$ consists of $x,y$ and all vertices
lying in some $x-y$ geodesic of $G$. For $S\subseteq V(G)$, the
\emph{geodetic closure} $I[S]$ of $S$ is the union of all  closed
intervals $I[u,v]$ over all pairs $u,v\in S$, i.e.,
$I[S]=\bigcup_{u,v\in S}I[u,v]$. The set $S$  is called
\emph{geodetic } if $I[S] = V(G)$ and it is said to be \emph{convex}
if $I[S] = S$. The \emph{convex hull} $CH(S)$ of $S$ is the smallest convex set containing $S$. If we define $I^0[S]=S, I^i[S]=I[I^{i-1}[S]$ for every $i\ge1$, then $CH(S)=I^r[S]$, for some $r\ge0$. 
The set $S$ is said to be a \emph{hull
set} if its convex hull $CH(A)$ is the whole vertex set
$V(G)$. The \emph{geodetic number} $g(G)$ and the \emph{hull
number} $h(G)$ are the minimum cardinality of a geodetic set and a
hull set, respectively \cite{ES:85,HLT:93}. Certainly, every
geodetic set is a hull set, and hence, $h(G)\leq g(G)$. In  Table \ref{t1}, both the geodetic number and the hull number of some families of graphs are shown.

\begin{table}[htbp]
\begin{center}\small
\begin{tabular}{||c||c|c|c|c|c|c|c|c||}
\hline
$G$&$P_n$&$C_{2l}$&$C_{2l+1}$&$T_n^h$&$K_n$ & $K_{p,n-p}$ & $S_{1,n-1}$  & $W_{1,n-1}$ \\
\hline\hline

$h(G)$&$2$&$2$&$3$& $h$ & $n$ & $2$ & $n-1$ & $\lceil\frac{n-1}{2}\rceil$\\
\hline

$g(G)$&$2$&$2$&$3$& $h$ &$n$&$\min\{4,p\}$  & $n-1$   & $\lceil\frac{n-1}{2}\rceil$\\
\hline
\end{tabular}
\end{center}
\vspace{-.4cm}\caption{Hull number and geodetic number of some graph classes.\label{t1}}
\end{table}

\begin{rem}\label{notaciones} In the rest of this paper,  $P_n$, $C_n$ and $K_n$ denote the path, cycle and complete graph of order $n$, respectively. In all cases, unless otherwise stated, the set of vertices is $\{0,1,\cdots,n-1\}$. In addition, $K_{p,n-p}$, $S_{1,n-1}$, $W_{1,n-1}$  denote the complete bipartite graph (being its smallest stable set of order $p\ge 2$), star and wheel of order $n$, whereas $T_n^h$ represents an arbitrary tree of order $n$ with $h$ leaves. Finally, in the sequel, $G$ and $H$ denote a pair of nontrivial connected graphs.
\end{rem}

 A vertex $v\in V(G)$ is a \emph{simplicial vertex} if the subgraph induced by its
neighborhood $N(v)=\{u:uv\in E(G)\}$ is a complete graph. It is easily seen that every hull
set, and  hence every geodetic set, must contain the set $Ext(G)$ of
simplicial vertices of $G$. A graph $G$ is called {\it extreme
geodesic} if the set of its simplicial vertices is geodetic (see \cite{CZ:04}).  Note that,
in this case, (1) the set $Ext(G)$ is the unique minimum geodetic
set (and also the unique minimum hull set) and (2)
$h(G)=g(G)=|Ext(G)|$. Trees and complete graphs are basic examples
of  extreme geodesic graphs.

\section{Strong product of graphs: general results}
\label{gs}

The \emph{strong product}  of graphs $G$ and $H$, denoted by $G \boxtimes
H$,  is the graph with the vertex set $V(G)\times V(H)=\{(g,h): g\in V(G), h\in V(H)\}$ in which vertices  $(g,h)$ and $(g',h')$ are adjacent whenever (1) $g = g'$ and $hh' \in E(H)$, or (2) $h = h'$ and $gg'
\in E(G)$, or (3) $gg' \in E(G)$ and $hh' \in E(H)$. 

The most important metric property of the strong product operation, relating the distance between two arbitrary vertices of an strong product graph to the distances between the corresponding projections in its factors,  is shown next.

\begin{lem}{\rm(\cite{IK:00})}\hspace{.1cm}
If $(g,h),(g',h')\in V(G \boxtimes
H)$, then $d_{G \boxtimes H}
((g,h), (g',h')) = max \{d_G (g, g'), d_H (h, h')\}$. 
Hence,  $diam(G \boxtimes H)=\max\{diam(G),diam(H)\}$.
\end{lem}

In this section, we firstly present some lemmas in order to show the behavior of the closed interval operator with respect to the strong graph operation, and next, we analyze in which way, both geodetic and hull sets of the strong product  of two graphs, are related to geodetic and hull sets of each factor, in both directions.

In the sequel, $p_G(S)$ and $p_H(S)$ denote the projections of a set of vertices  $S\subseteq V(G\boxtimes H)$ onto $G$ and $H$, respectively.

\begin{lem}
\label{projcamino}
Let $u=(g,h),v=(g',h')\in V(G \boxtimes
H)$ such that $d_{G \boxtimes H}
(u,v) = d_G (g, g')=l$. If $\gamma$ is a  $(g,h)-(g',h')$ geodesic, then the projection of $\gamma$ onto $G$ is a  $g-g'$ geodesic of length $l$.
\end{lem}
\begin{proof} 
If $V(\gamma)=\{(g,h),(g_1,h_1),\ldots,(g_{l-1},h_{l-1}),((g',h')\}$, then its projection  into $G$ is $p_G(V(\gamma))=\{g,g_1,\ldots,g_{l-1},g'\}$. Since $d_{G \boxtimes H}
((g,h), (g',h')) = d_G (g, g')$, $p_G(V(\gamma))$ does not contain repeated vertices, which means that every pair of consecutive vertices are adjacent, i.e., $p_G(V(\gamma))$ is  the vertex set of a $g-g'$ geodesic in $G$.
\end{proof}

\begin{lem}\label{lem.projmax}
  Let  $u=(g_1,h_1),v=(g_2,h_2)\in V(G\boxtimes H)$ such that $d_{G \boxtimes H}
(u,v) = d_G (g_1, g_2)=l$.
  Then, $$I[u,v]=\{ (g,h) : g\in I[g_1,g_2], d_H(h_1,h)\le d_G(g_1,g), d_H(h,h_2)\le d_G(g,g_2) \}.$$
\end{lem}
\begin{proof}
  Let $w=(g,h)$ be a vertex belonging to $I[u,v]$. By Lemma~\ref{projcamino}, the projection  of every $u-v$ geodesic onto $G$ is a $g_1-g_2$ geodesic, which means that $g\in I[g_1,g_2]$. If  $d_G(g_1,g)> d_H(h_1,h)$, then
\begin{align*}
    d_G(g_1,g_2)&=d_{G \boxtimes H}(u,v)=d_{G \boxtimes H}(u,w)+d_{G \boxtimes H}(w,v)\\&=d_H(h_1,h)+\max \{ d_G(g,g_2), d_H(h,h_2) \} \\& > d_G(g_1,g)+d_G(g,g_2)=d_G(g_1,g_2),
  \end{align*}
 which is a contradiction. Similarly,  a contradiction is obtained  by assuming that $d_H(h,h_2)> d_G(g,g_2)$.
  
 Conversely, suppose  that $w=(g,h)$ is a vertex belonging to $V(G\boxtimes H)$ such that $g\in I[g_1,g_2]$, 
 $r=d_H(h_1,h) \le d_G(g_1,g)=k$, $s=d_H(h,h_2) \le  d_G(g,g_2)=l-k$. Let $\rho$ be a $g_1-g_2$ geodesic passing through $g$ such that $V(\rho)=\{z_0,z_1,\ldots,z_l\}$, $z_0=g_1$, $z_k=g$  and $z_l=g_2$. Let $\mu_1$ be a $h_1-h$ geodesic such that $V(\mu)=\{x_0,x_1,\ldots,x_r\}$, $x_0=h_1$ and $x_r=h$. Let $\mu_2$ be a $h-h_2$ geodesic such that $V(\mu_2)=\{y_0,y_1,\ldots,y_s\}$, $y_0=h$ and $y_s=h_2$. It is straightforward to check that 
 $$\{ (z_0,x_0),(z_1,x_1),\dots ,(z_r,x_r),\dots , (z_k,x_r),(z_{k+1},y_1), \dots ,(z_{k+s},y_s),\dots  ,(z_l,y_s)  \} $$
is the vertex set of a $u-v$ geodesic passing through $w$, which means that $w\in I[u,v]$.
\end{proof}

\begin{figure}[htb]
  \begin{center}
  \includegraphics[width=0.6\textwidth]{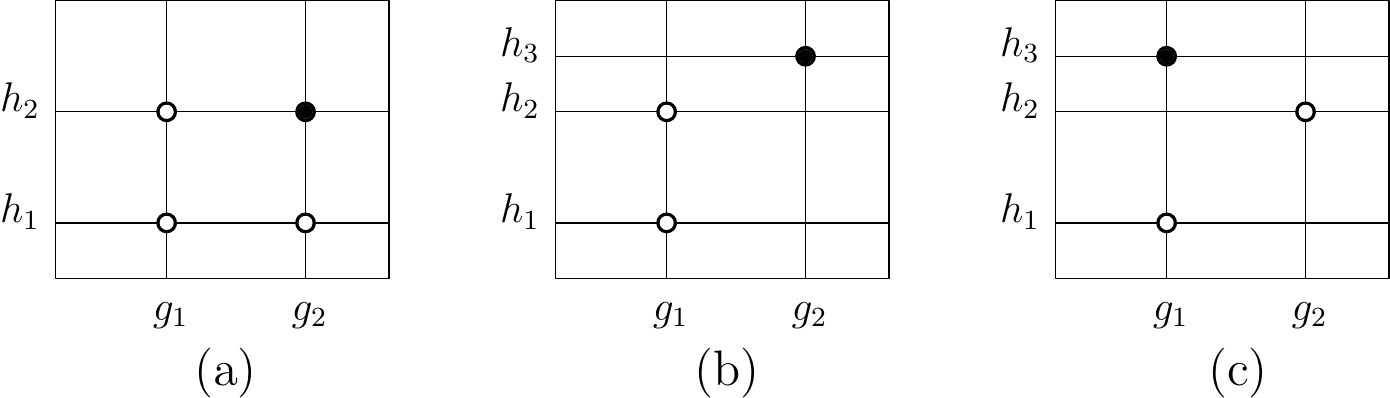}
  \caption{In each figure, the black vertex does not belong to the geodetic closure of white vertices}
  \label{fig.cor}
  \end{center}
\end{figure}

\begin{lem}
\label{lemazo} Let  $S_1\times S_2\subseteq V(G \boxtimes H)$ a set of vertices of cardinality 6, where $S_1=\{g_1,g_2\}\subseteq V(G)$ and $S_2=\{h_1,h_2,h_3\}\subseteq V(H)$. Then

\begin{itemize}

\item[(i)] $(g_2,h_2)\not\in I[(g_1,h_1),(g_1,h_2),(g_2,h_1)]$

\item[(ii)] If $h_3\not\in I[h_1,h_2]$, then $(g_2,h_3)\notin I[(g_1,h_1),(g_1,h_2)]$. 

\item[(iii)] If $h_3\not\in I[h_1,h_2]$, then $(g_1,h_3)\notin I[(g_1,h_1),(g_2,h_2)]$. 

\end{itemize}
\end{lem}
\begin{proof} 

(i) Observe that $d((g_1,h_1),(g_1,h_2))=d(h_1,h_2)$. Hence, according to 
Lemma~\ref{projcamino},  every $(g_1,h_1)-(g_1,h_2)$ geodesic may not pass through $(g_2,h_2)$.  Similarly, it  is  proved that $(g_2,h_2)\not\in I[(g_1,h_1),(g_2,h_1)]$ and $(g_2,h_2)\not\in I[(g_1,h_2),(g_2,h_1)]$ (see Figure~\ref{fig.cor}(a)).

(ii) Observe that $d((g_1,h_1),(g_1,h_2))=d(h_1,h_2)$. Hence, according to 
Lemma~\ref{projcamino}, the projection  onto $H$ of a $(g_1,h_1)-(g_1,h_2)$ geodesic passing through $(g_2,h_3)$ is a $h_1-h_2$ geodesic passing through $h_3$, contradicting the hypothesis $h_3\not\in I[h_1,h_2]$  
(see Figure~\ref{fig.cor}(b)).

(iii) Suppose that  $(g_1,h_3)\in I[(g_1,h_1),(g_2,h_2)]$. 

If $d((g_1,h_1),(g_2,h_2))=d(g_1,g_2)$,
then, according to 
Lemma~\ref{projcamino}, every $(g_1,h_1)-(g_2,h_2)$ geodesic may not pass through $(g_1,h_3)$.

If $d((g_1,h_1),(g_2,h_2))=d(h_1,h_2)$, 
then the projection onto $H$ of a $(g_1,h_1)-(g_2,h_2)$ geodesic passing through $(g_1,h_3)$ is a $h_1-h_2$ geodesic passing through  $h_3$, which contradicts the hypothesis $h_3\notin I[h_1,h_2]$ (see Figure~\ref{fig.cor}(c)).
\end{proof}

\begin{lem}
\label{lem.intervalos}
Let  $S_1\subseteq V(G)$ and $S_2\subseteq
V(H)$. Then, 
for every  integer $r\geq 1$, $I^r[S_1]\times I^r[S_2]\subseteq I^r[S_1\times S_2]$.
\end{lem}
\begin{proof}
We proceed by induction on $r$. Suppose that $r=1$ and take a vertex $(g,h)\in I[S_1]\times I[S_2]$. Since $g\in I[S_1]$, then
$g\in I[g',g'']$ for some $g',g''\in V(S_1)$, and thus
$d(g',g'')=d(g',g)+d(g,g'')$. Similarly,
$d(h',h'')=d(h',h)+d(h,h'')$ for some $h',h''\in V(S_2)$. We may assume without loss of generality that $d(g',g)\leq d(g,g'')$, $d(h',h)\leq d(h,h'')$ and  $d(g',g)\leq d(h',h)$. 
Then, $d((g',h'),(g,h))=d(h',h)$ and  $d((g,h),(g',h''))=d(h,h'')$, which means that
$$ d((g',h'),(g',h''))=d(h',h'')=d(h',h)+d(h,h'')=d((g',h'),(g,h))+d((g,h),(g',h''))$$
In other words, $(g,h)\in I[(g',h'),(g',h'')]\subseteq I[S_1\times S_2]$.

Assume then that $r>1$. By the inductive hypothesis, $I^{r-1}[S_1]\times I^{r-1}[S_2]\subseteq I^{r-1}[S_1\times S_2]$. Hence,
$ I^r[S_1]\times I^r[S_2]=I[I^{r-1}[S_1]]\times I[I^{r-1}[S_2]]\subseteq I[I^{r-1}[S_1]\times I^{r-1}[S_2]]\subseteq I[I^{r-1}[S_1\times S_2]]=I^r[S_1\times S_2]$. 
\end{proof}

As a direct consequence of this lemma, the following result is obtained.

\begin{prop}
\label{prop:proposicion1} Let  $S_1\subseteq V(G)$ and $S_2\subseteq
V(H)$. If  $S_1$ is a  geodetic (resp. hull) set of  $G$ and $S_2$ is a geodetic  (resp. hull) set of $H$, then   $S_1\times
S_2$ is a geodetic (resp. hull) set of  $G\boxtimes H$.
\end{prop}

\begin{proof}
Let $r,s$ be positive integers such that $I^r[S_1]=V(G)$ and $I^s[S_2]=V(H)$. We may suppose wlog that $r\leq s$. Then, $V(G\boxtimes H)=V(G)\times V(H)=I^s[S_1]\times I^s[S_2]\subseteq  I^s[S_1\times S_2]$.
\end{proof}

\begin{prop}
\label{prop:proposicion2}
Let $S\subseteq V(G\boxtimes H)$ be a geodetic set of  $G\boxtimes H$.  Then,  either the projection  of $S$ onto $G$ or the projection  of $S$ onto $H$ is geodetic.
\end{prop}

\begin{proof}
Assume that neither
$S_1=p_G(S)$ nor $S_2=p_H(S)$ is geodetic and consider $g\in V(G)\backslash I[S_1]$ and
$h\in V(H)\backslash I[S_2]$. As $(g,h)\in
I[S]=V(G\boxtimes H)$, then  $(g,h)\in I[(g',h'),(g'',h'')]$ for some $(g',h'),(g'',h'')\in S$. Hence,
$ d((g',h'),(g'',h''))=d((g',h'),(g,h))+d((g,h),(g'',h''))$.

On the other hand, as $g\notin I[g',g'']$ and $h\notin I[h',h'']$, we have that
$d(g',g'')<d(g',g)+d(g,g'')$ and
$d(h',h'')<d(h',h)+d(h,h'')$. Hence,
$$\max \{d(g',g''),d(h',h'')\}<\max \{d(g',g)+d(g,g''),d(h',h)+d(h,h'')\}\leq $$
$$\leq \max \{d(g',g),d(h',h)\}+\max\{d(g,g''),d(h,h'')\}=d((g',h'),(g,h))+d((g,h),(g'',h''))$$
which contradicts the previous expression for the distance between $(g',h')$ and $(g'',h'')$.
\end{proof}

This property is far from being true for hull sets, as it is shown
in the next example.

\begin{exa}\label{c5c7}
It is straightforward to prove that if $V(C_5)=\{ u_1,\dots , u_5\}$ and $V(C_7)=\{ v_1,\dots , v_7\}$, then (1) $\{(u_1,u_2\}$ is not a hull set of $C_5$, (2) $\{(v_1,v_4\}$ is not a hull set of $C_7$, and (3) $\{(u_1,v_1), (u_2,v_4)\}$ is a  hull set of $C_5\boxtimes C_7$.
\end{exa}

\section{Geodetic and hull numbers: bounds}

In this section, we study the behavior of  both the geodetic and the hull numbers
with respect to  the  strong product operation for
graphs, in terms of its factors. More precisely, we obtain bounds, 
and we give some examples showing that all of them
are sharp.

\begin{lem}\label{2de3}  Let $\{h_1,h_2,h_3\}$ a 3-vertex set of a graph $H$.
If $h_1\in I[h_2,h_3]$, then $h_2\notin I[h_1,h_3]$ and $h_3\notin I[h_1,h_2]$.
\end{lem}
\begin{proof}
Assume on the contrary that, for example,  $h_2\in I[h_1,h_3]$. Then, if $d(h_1,h_2)=x$, $d(h_1,h_3)=y$ and $d(h_2,h_3)=z$, we have that $x+y=z$ and $x+z=y$, i.e.,  $d(h_1,h_2)=0$, a contradiction.
\end{proof}

\begin{prop}\label{prop.gno3}
Let $G$ and $H$ be nontrivial graphs. Then, $g(G\boxtimes H)\geq 4$.
\end{prop}

\begin{figure}[htb]
  \begin{center}
  \includegraphics[width=0.7\textwidth]{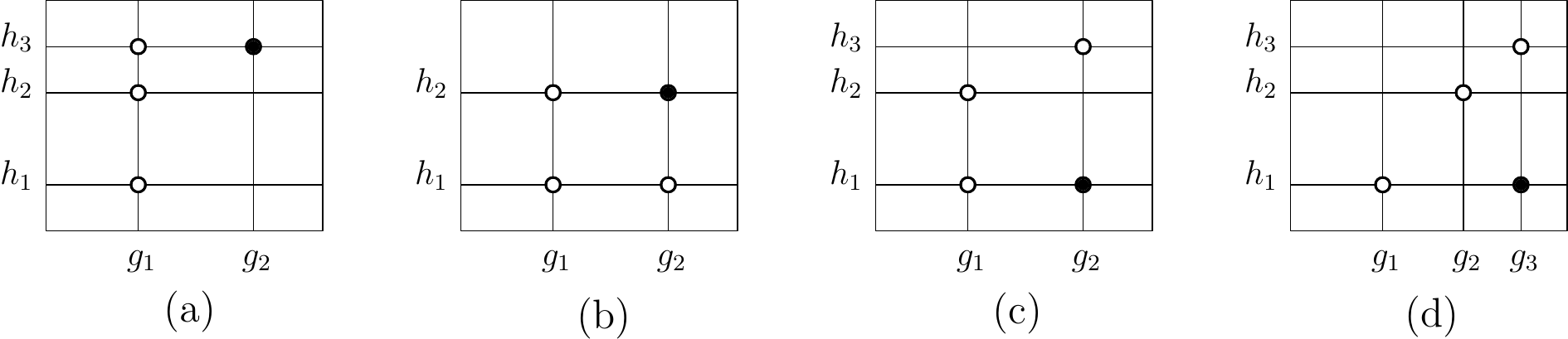}
  \caption{In each figure, the black vertex does not belong to the geodetic closure of white vertices}
  \label{4casos}
  \end{center}
\end{figure}

\begin{proof}

Let us see that every subset  $S$ of $V(G\boxtimes H)$ having at most 3 vertices is not geodetic. Suppose on the contrary that $S$ is a geodetic set of cardinality 3. Without loss of generality, we may assume that $|p_G(S)|\leq|p_H(S)|$. We consider different cases.

{\bf Case 1.} $|p_G(S)|=1$: In other words, $S=\{ (g_1,h_1), (g_1,h_2), (g_1,h_3)\}$ and  $|p_H(S)|=3$.  According to 
Lemma~\ref{2de3}, we may assume w.o.l.g. that $h_3\notin I[h_1,h_2]$  and from Lemma~\ref{lemazo}(i,ii), we derive that $(g_2,h_3)\notin I[S]$ for any vertex $g_2\not= g_1$  (see Figure \ref{4casos}(a)).

{\bf Case 2.} $|p_G(S)|=|p_H(S)|=2$: In other words, $S=\{ (g_1,h_1), (g_1,h_2), (g_2,h_1)\}$, being $g_1\neq g_2$,  and $h_1\neq h_2$. From Lemma~\ref{lemazo}(i), we derive that   $(g_2,h_2)\notin I[S]$  (see Figure \ref{4casos}(b)).

{\bf Case 3.} $|p_G(S)|=2$ and $|p_H(S)|=3$:  In other words, $S=\{ (g_1,h_1), (g_1,h_2), (g_2,h_3)\}$, being $g_1\neq g_2$,  and $h_1,h_2,h_3$ three  diferent vertices  of $H$. According to 
Lemma~\ref{2de3}, we may assume w.o.l.g. that $h_1\notin I[h_2,h_3]$. From Lemma~\ref{lemazo}(i,iii), we  derive that $(g_2,h_1)\notin I[S]$  (see Figure \ref{4casos}(c)).

{\bf Case 4.} $|p_G(S)|=|p_H(S)|=3$: In other words, $S=\{ (g_1,h_1), (g_2,h_2), (g_3,h_3)\}$, being $g_1, g_2,g_3 $ three different vertices of $G$ and $h_1,h_2,h_3$ three different vertices of $H$. According to 
Lemma~\ref{2de3}, we may assume w.o.l.g. that $h_1\notin I[h_2,h_3]$ and $g_3\notin I[g_1,g_2]$. From Lemma~\ref{lemazo}(i,iii), we  derive that  $(g_3,h_1)\notin I[S]$ (see Figure \ref{4casos}(d)).
\end{proof}

As a direct consequence of Propositions~\ref{prop:proposicion1} and ~\ref{prop:proposicion2}, we derive bounds for
the geodetic number of the strong product of two graphs , in terms
of the geodetic numbers of its factor graphs.

\begin{thm}
\label{teo.geodetico}
For any two graphs $G$ and $H$,
$ \min \{g(G),g(H)\}\leq g(G\boxtimes H)\leq g(G)g(H).$

\noindent Furthermore, both bounds are sharp.
\end{thm}
\begin{proof}
First, we prove the upper bound. Let $S_1$ and $S_2$ be
geodetic sets of $G$ and $H$ with minimum cardinality, that is, such that
$|S_1|=g(G)$ and $|S_2|=g(H)$. By
Proposition~\ref{prop:proposicion1}, $S_1\times S_2$ is a geodetic
set of $G\boxtimes H$ with cardinality $|S_1\times
S_2|=|S_1||S_2|=g(S_1)g(S_2)$. Hence, $g(G\boxtimes H)\leq g(G)g(H)$.

To prove  the lower bound, take a minimum geodetic set $S$ of
$G\boxtimes H$. According to Proposition~\ref{prop:proposicion2}, we
may suppose, without loss of generality, that $p_G(S)$  is a
geodetic set of $G$. Hence: $\min \{g(G),g(H)\}\leq g(G) \leq
|p_G(S)| \leq |S| = g(G\boxtimes H)$.

To show the sharpness of the upper bound, take $G=K_m$ and
$H=K_n$. Then, $g(K_m\boxtimes K_n)=g(K_{mn})=mn=g(K_m)g(K_n)$.
Finally, to show the sharpness of the lower bound, take $G=K_{r,s}$  a complete bipartite graph and
$H=K_n$, with $r,s,n\geq 4$. Then, as it will be shown in the next section (see Example~\ref{bicliques}), $G(K_{r,s}\boxtimes K_n)=4=\min\{ g(K_{r,s}), g(K_n)\}$.
\end{proof}

\begin{thm}
\label{teo.hull}
For any  two nontrivial graphs $G$ and $H$,  $2\le h(G\boxtimes H)\leq h(G)h(H)$. Furthermore, both bounds are sharp.
\end{thm}
\begin{proof}
First, we prove the upper bound. Let $S_1$ and $S_2$ be
hull sets of $G$ and $H$ with minimum cardinality, that is, such that
$|S_1|=h(G)$ and $|S_2|=h(H)$. By
Proposition~\ref{prop:proposicion1}, $S_1\times S_2$ is a hull
set of $G\boxtimes H$ with cardinality $|S_1\times
S_2|=|S_1||S_2|=h(S_1)h(S_2)$. Hence, $h(G\boxtimes H)\leq h(G)h(H)$.

To prove the sharpness of this  bound, take $G=K_m$ and $H=K_n$ and notice that $h(K_m\boxtimes K_n)=h(K_{mn})=mn=h(K_m)h(K_n)$

Finally, the lower bound is a directe consequence of the fact that $h(G)=1$ if and only if $G=K_1$. As for its sharpness, it is straightforward to check that $\{(0,0),(0,2)\}$ is a hull set of $P_2\boxtimes C_4$.
\end{proof}

\begin{rem} Conversely to the geodetic case,  the claim $\min \{h(G),h(H)\}\leq h(G\boxtimes H)$ is far from being true in general. A simple counterexample is shown in Example~\ref{c5c7}.
\end{rem}

\begin{lem}
\label{lem:rodaja}
Let  $G$ and $H$ be two  graphs such that $Ext(G)=\emptyset$.
If $S$ is a hull set of $G$ and $x$ is an arbitrary vertex of $H$, then  $S\times \{ x \}$ is a hull set of $G\boxtimes H$.
\end{lem}
\begin{proof}
We prove by induction on $m\ge 0$ that for every vertex $h\in V(H)$ such that $d(x,h)=m\ge 0$, if $g\in V(G)$, then the   vertex $(g,h)$ is in the convex hull of $S\times \{ x\}$.

For $m=0$, the condition $d(x,h)=m\ge 0$ implies $h=x$. Since  $S$ is a hull set of $G$, for every $g\in V(G)$ we have $g\in  I^r[S]$ for some $r\ge 0$. By lemma~\ref{lem.intervalos}, $(g,x)\in I^r[S]\times I^r[\{ x\}]=I^r[S\times \{ x\} ]$, and consequently, $(g,x)$ is in the convex hull of $S\times \{ x \}$.

Suppose now $m>0$ and consider a vertex $h\in V(H)$ with $d(x,h)=m>0$. Take a vertex $h'\in V(H)$  such that $d(h,h')=1$ and  $d(h',x)=m-1$. Since $G$ has no simplicial vertices, for every vertex $g\in V(G)$ there exist vertices $g_1,g_2$ in $G$ adjacent to $g$ such that $d(g_1,g_2)=2$. Thus $d((g_1,h'),(g_2,h'))=2$, $d((g,h),(g_1,h'))=1$ and $d((g,h),(g_2,h')=1$, that is,  $(g,h)$ is in a geodesic
between $(g_1,h')$ and $(g_2,h')$.  By inductive hypothesis, $(g_1,h')$ and $(g_2,h')$ are in the convex hull of $S\times \{x\}$. Therefore, $(g,h)$ is in the convex hull of $S\times \{ x \}$.
\end{proof}

As a consequence of the preceding lemma we obtain the following upper bound for the hull number of the strong product of two graphs, if at least one of them has no simplicial vertices. 

\begin{thm}
\label{teo.hull2}
Let  $G$ and $H$ be two graphs such that $Ext(G)=\emptyset$. Then, $h(G\boxtimes H)\le h(G)$
\end{thm}

Certainly, this last bound is also sharp. Consider, for example the strong product graph $C_m\boxtimes C_n$, being  both $m$ and $n$ even. As it will be shown in the next section (see Proposition~\ref{toruses}), $h(C_{m}\boxtimes C_n)=2=h(C_{m})$.

\section{Exact values}

In this section, we approach the calculation of   the geodetic  and the hull numbers of a some  strong product graphs, where at least one of the factors is either a complete graph or a cycle or a path. 
We begin by showing a result involving extreme geodesic graphs.


\begin{prop}\label{extreme0}
Two graphs  $G$ and $H$ are extreme geodesic  if and only if $G\boxtimes H$ is an extreme geodetic graph.
\end{prop}
\begin{proof}
Observe  that a
vertex $(g,h)$ is a simplicial vertex of  $G\boxtimes H$ if and only
if both $g$ and $h$ are simplicial vertices of $G$ and $H$,
respectively, i.e., $Ext(G\boxtimes H)=Ext(G)\times Ext(H)$. As a
direct consequence of this equality and
Proposition~\ref{prop:proposicion1}, we  have that  two  graphs $G$
and $H$ are extreme geodesic if and only if $G\boxtimes H$ is
extreme geodesic.
\end{proof}

\begin{cor}\label{extreme}
If both $G$ and $H$ are extreme geodesic graphs, then $h(G\boxtimes H)=g(G\boxtimes
H)=g(G)g(H)=h(G)h(H).$
\end{cor}

As a direct consequence of Corollary~\ref{extreme},
the results shown in Table~\ref{t3} are obtained.

\begin{table}[htbp]
\begin{center}\small
\begin{tabular}{||c||c|c|c||}
\hline

$G/H$ & $P_n$ & $T_n^k$ & $K_n$ \\
\hline\hline

$P_m$ & 4 & $2k$ & $2n$ \\
\hline

$T_m^h$ & $2h$ & $hk$ & $hn$ \\
\hline

$K_m$ & $2m$ & $mk$ & $mn$ \\
\hline\hline

\end{tabular}
\end{center}
\vspace{-.4cm}\caption{ Hull  and geodetic numbers of the strong product of some extreme geodesic  graphs.\label{t3}}
\end{table}

Certainly, cycles are graphs without simplical vertices, and hence they are not extreme geodesic graphs. This means that the calculation of the geodetic and the hull numbers of strong product graphs of the form $G\boxtimes C_n$, requires a different approach to the previous one.  The rest of this section is devoted to this issue.

\begin{defi}
Let be  $S$ a set of vertices in a graph $G$. Then, $S$ is said to satisfy condition
\begin{itemize}
\item[(A)] if, for every vertex $x\in S$,  there exist two vertices $y,z\in S-x$ such that $x\in I[y,z]$.
\item[(B)] if there are two vertices $x,y\in S$ such that $x\notin
I[S-x]$ and $y\notin I[S-y]$.
\end{itemize}
\end{defi}

\begin{lem}\label{lemaSs}
Let G be a graph having a geodetic set $S$ satisfying the condition (A).
Then, for every vertex $k\in V(K_n)$, $S\times \{k\}$ is a geodetic set of $g(G~\boxtimes~K_n)$.
\end{lem}
\begin{proof}
Take an arbitrary vertex $(g,h)\in V(G\boxtimes K_n)$. This means that there exists a pair of vertices $s,s' \in S\setminus \{g\}$ such that $g\in I[s,s']$. Hence, 
$d((s,k), (s',k))=d(s,s')=d(s,g)+d(g,s')=d((s,k),(g,h))+d((g,h),(s',k))$, i.e., 
$(g,h)\in I[(s,k),(s',k)]\subset I[S]$, as desired.
\end{proof}

\begin{prop}
Let G be a graph with a minimum geodetic set $S$ satisfying condition (A).
Then, for every positive integer $n$, $g(G~\boxtimes~K_n)~=~g(G)$.
\end{prop}
\begin{proof}
As a corollary of Lemma~\ref{lemaSs} we have that $g(G\boxtimes K_n) \leq g(G)$. To get the equality,
suppose  that there exists a geodetic set $R=\{ (g_1, k_1), (g_2, k_2),\dots (g_m,k_m)\}$ in $G\boxtimes K_n$ such that $m=|R| <|S|=g(G)$. Consider the set $R'=\{ (g_1, k_1), (g_2, k_1),\dots (g_m,k_1)\}$. For every vertex  
$(g,k)\in G~\boxtimes~K_n$ we have that $(g,h)\in I[(g_i, k_i), (g_j, k_j)]$ for some $i,j\in\{1,\ldots,m\}$. Hence, $g\not=g_i\not=g_j\not=g$ and
$d((g_i, k_1), (g_j, k_1))=d(g_i,g_j)=d((g_i, k_i), (g_j, k_j))=
d((g_i, k_i), (g, k))+d((g, k), (g_j, k_j))=d(g_i,g)+d(g,g_j)=d((g_i, k_1), (g, k))+d((g, k), (g_j, k_1))$. In other words, $(g,k)\in I[(g_i, k_1),(g_j, k_1)]\subseteq I[R']$. We have thus proved that $R'$ is also a geodetic set of $G~\boxtimes~K_n$. 
Furthermore, as a direct consequence of  Proposition~\ref{prop:proposicion2}, we conclude that the projection $p_G(R')$ is a geodetic set of $G$, from which it follows that  $|p_G(R')|= |R'|\le|R|<|S|=g(G)$, a contradiction.
\end{proof}

\begin{exa}\label{bicliques} Consider the complete bipartite graph $K_{r,s}$, whith $2\le r\le s$. Notice that if $V(K_{r,s})=\{u_1,\ldots,u_r\}\cup\{v_1,\ldots,v_s\}$, then the set $\{u_1,u_2,v_1,v_2\}$ is a minimum geodetic set satisfying condition (A). Hence,
$g(K_{r,s}~\boxtimes~K_n)~=~g(K_{r,s})=4$.
\end{exa}

\begin{prop}
Let $n\ge4$ be an even  integer and let $G$ be a graph of order $m\ge2$. If $G$ is either a path $P_m$ or a complete graph $K_m$, then: 
$g(G\boxtimes C_n)=4$ and $h(G\boxtimes C_n)=2$.
\end{prop}
\begin{proof}
The equality $h(G\boxtimes C_n)=2$ is a direct consequence of Theorem~\ref{teo.hull2}. The equality $g(P_m\boxtimes C_n)=4$ is a corollary of Proposition~\ref{prop.gno3} and the upper bound shown in Theorem~\ref{teo.geodetico}. Finally, to prove that $g(K_m\boxtimes C_n)=4$  it is enough to consider again Proposition~\ref{prop.gno3} and to notice that the set $S=\{ 0, 1, \frac{n}{2}, \frac{n+2}{2}\}$ is a (not minimum) geodetic set in $C_n$ satisfying condition (A).
\end{proof}

\begin{lem}\label{vertincycles}  
Let $h\ge 2$ be an integer and let $S$ be a set of vertices in the cycle $C_{2h+1}$.  If $2\le|S|\le4$, then $S$
satisfies condition (B).
\end{lem}

\begin{proof}
Certainly, this statement  is trivial if $|S|=2$. For $|S|=3$,  we may assume that
$S=\{0,i,j\}$, where $0< i< j\le 2h$. If $i> h$, then $0\notin
I[i,j]$ and $i\notin I[0,j]$. If $j\le h$, then $0\notin I[i,j]$
and $j\notin I[0,i]$. If $i\le h$ and $j\ge h+1$, then $i\notin
I[0,j]$ and $j\notin I[0,i]$.

\begin{figure}[h]
\begin{center}
  \includegraphics[width=0.7\textwidth]{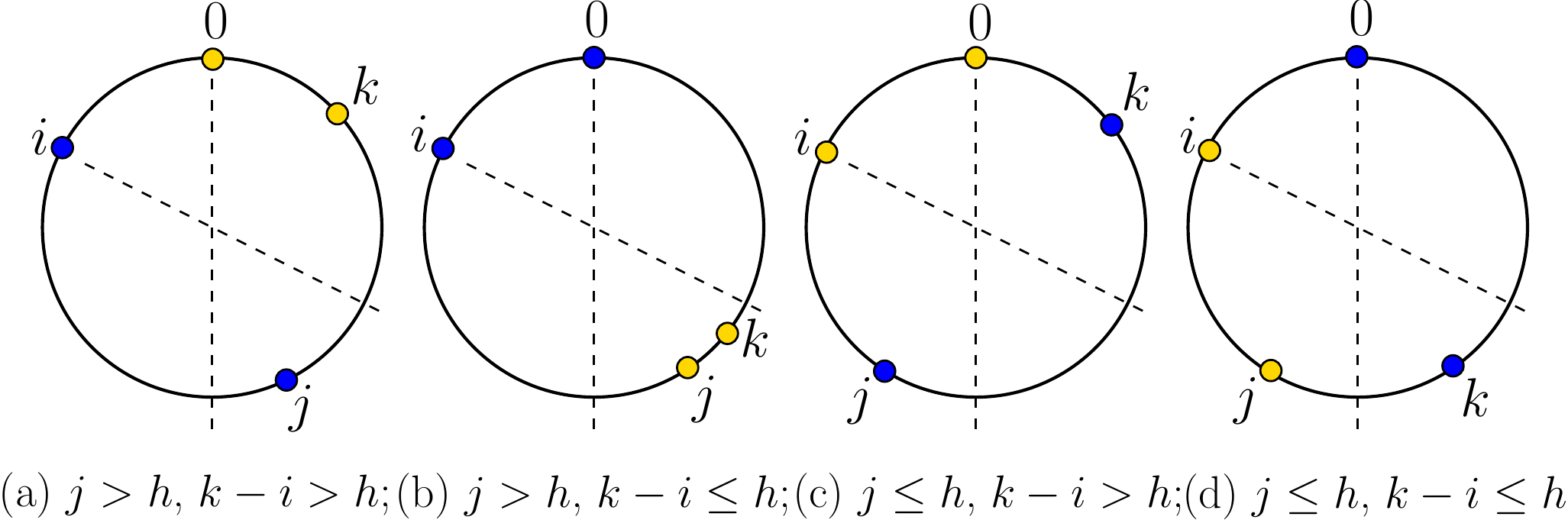}
 \caption{In all cases, $S=\{0,i,j,k\}$ and $0< i< j<k
\le 2h$.}
  \label{rodones}
\end{center}
\end{figure}

For $|S|=4$,  we may assume that $S=\{0,i,j,k\}$, where $0< i< j<k
\le 2h$. If $j>h$ and $k-i>h$, then $i\notin I[S-i]$ and $j\notin
I[S-j] $. If $j>h$ and $k-i\le h$, then $i\notin I[S-i]$ and
$0\notin I[S-0] $. If $j\le h$ and $k-i>h$, then $k\notin I[S-k]$
and $j\notin I[S-j] $. If $j\le h$ and $k-i\le h$, then $k\notin
I[S-k]$ and $0\notin I[S-0]$ (see Figure \ref{rodones}).
\end{proof}

\begin{prop}\label{ultima} Let $G$ be a nontrivial graph such that 
every set of vertices $S\subseteq V(G)$ of cardinality $2\le|S|\le4$ satisfies condition (B). Then, for every integer $k\ge 2$, $g(G \boxtimes C_{2k+1})\geq 5$.
\end{prop}
\begin{proof} Denote $H=C_{2k+1}$ and assume  that $g(G \boxtimes H)= 4$, i.e., that $S$ is a geodetic set  of cardinality 4. 
Observe that $1\le |p_G(S)|\le 4$ and $1\le |p_H(S)|\le 4$. We distinguish three cases.

{\bf Case 1.} $|p_G(S)|=1$ or $|p_H(S)|=1$: 
If $|p_G(S)|=1$, as $H$ satisfies condition (B), there exists a vertex $h\in
p_H(S)$ s.t. $h\notin I[p_H(S)-h]$. If $p_G(S)=\{g_1\}$ and $g_1g\in E(G)$ then,  according to Lemma \ref{lemazo}, $(g,h)\notin I[S]$ (see Figure \ref{puntitos}(a)). The case $|p_H(S)|=1$ is similarly proved (see Figure \ref{puntitos}(b)).
\begin{figure}[h]
\begin{center}
  \includegraphics[width=0.8\textwidth]{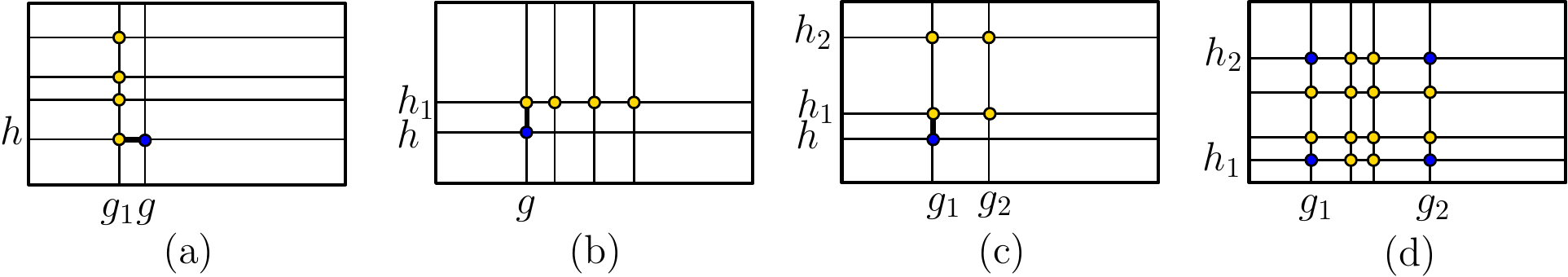}
 \caption{In all cases, each dark vertex is not in the geodetic closure of the remaining vertices.}
  \label{puntitos}
\end{center}
\end{figure}

{\bf Case 2.} $|p_G(S)|= |p_H(S)|=2$: If $p_G(S)=\{ g_1,g_2\}$ and
$p_H(S)=\{ h_1,h_2\}$. If $h$ is the  vertex  adjacent to $h_1$ not belonging to
$I[h_1,h_2]$ then, by Lemma \ref{lemazo}, $(g_1,h)\notin I[S]$ (see Figure \ref{puntitos}(c)).

{\bf Case 3.} $|p_G(S)|\cdot|p_H(S)|>4$: 
As both $G$ and $H$ satisfy condition (B), 
there exist vertices $g_1,g_2$ in $G$ such that
$g_i\notin I[p_G(S)-g_i]$, for $i=1,2$, and vertices $h_1,h_2$ in $H$ such that $h_i\notin I[p_H(S)-h_i]$, for $i=1,2$. At
least one of the four vertices of  $\{ (g_i,h_j):i,j\in \{
1,2\} \}$, say $(g_1,h_1)$,  is not in $S$, as otherwise $|p_G(S)|\cdot|p_H(S)|=4$. 
Hence, by  Lemma \ref{lemazo}, $(g_1,h_1)\notin I[S]$ (see Figure \ref{puntitos}(d)). 
\end{proof}

\begin{prop} Let $n\ge5$ be an odd  integer. If $m\ge2$, then $g(K_m\boxtimes C_n)=5$ and $h(K_m\boxtimes C_n)=3$.
\end{prop}
\begin{proof}
Notice that  if  $n=2k+1$, then the set $S=\{ 0, 1, k, k+1, k+2\}$ is  a geodetic set of $C_n$ satisfying  condition (A), which means that $g(K_m\boxtimes C_n)\leq 5$.
The equality is directly derived from Proposition \ref{ultima}, since every set of vertices of $K_m$ trivially satisfies condition (B).

To prove that $h(K_m\boxtimes C_n)=3$ it suffices to see that $h(K_m\boxtimes C_n)>2$, as according to Theorem~\ref{teo.hull2}, $h(K_m\boxtimes C_n)\leq h(C_n)=3$. To this end, take an arbitrary set $R=\{ (i_1,j_1), (i_2,j_2)\}$ of cardinality  2 in $K_m\boxtimes C_n$. If $j_1=j_2$, then $CH(R)=R$, i.e., in this case $R$ is not a hull set of $h(K_m\boxtimes C_n)$. Assume thus that $j_1\not=j_2$, and wlog that $R=\{ (0,0), (i,h)\}$, where $i\in\{0,1\}$, $n=2k+1$  and $0<h\leq k$. Certainly, $CH(R)\setminus R=\bigcup_{j=1}^{h-1} (K_m\times\{j\})$, i.e.,  neither in this case $R$ is  a hull set of $h(K_m\boxtimes C_n)$.
\end{proof}

\begin{prop} For any $m,k\ge 2$, $5\le g(P_m\boxtimes C_{2k+1})\le 6$.
\end{prop}
\begin{proof} As an immediate consequence of Proposition~\ref{prop.gno3}   and Theorem~\ref{teo.geodetico} we obtain that $4\le g(P_m\boxtimes C_n)\le 6$. Moreover, observe that every set $S\subset V(P_m)$ such that $2\le|S|\le4$ satisfies condition (B), which according to Proposition \ref{ultima}, allows us to derive that  $g(P_m\boxtimes C_{2k+1})\geq 5$. 
\end{proof}

Let us remark that both bounds are sharp, since it is straightforward to check that $g(P_3\boxtimes C_{7})=5$ and $g(P_3\boxtimes C_{5})=6$.

\begin{prop}
For any $k,m\ge 2$, $h(C_{2k+1}\boxtimes P_m)=\left\{
                                              \begin{array}{ll}
                                                2, & \hbox{if }k\le m-2; \\
                                                3, & \hbox{if }k\ge m-1.
                                              \end{array}
                                            \right.$
\end{prop}
\begin{proof} Certainly, $2\le h(C_{2k+1}\boxtimes P_m)\le 3$, being the upper bound  a corollary of Theorem~\ref{teo.hull2}, whereas the lower bound is an immediate consequence of the fact that  $h(G)=1$ if and only if $G=K_1$.

By symmetry reasons, we label the vertex set of $C_{2k+1}$ as follows: $V(C_{2k+1})=\Lambda=\{ -k, \dots ,-1,0,$ $1,\dots ,k\}$, whereas the vertex set labeling is the usual one: $V(P_m)=\Pi=\{ 0,1,\dots ,m-1\}$.
In addition, we identify $V(C_{2k+1}\boxtimes P_m)$ with the grid $\mathcal{P}=\Lambda\times\Pi$ of points of the discrete plane $\mathbb{Z}^2$.

Suppose first that $k\le m-2$, $m$ is odd and 
take the vertices $u=(0,0)$ and $v=(0,m-1)$. Observe that, as shown in Figure \ref{fig.caminoHull2} (a), $I[u,v]$ is the set of points belonging to  the square determined by $\{u,v,x,y\}$, where $x=(-\frac{m-1}2,\frac{m-1}2)$ and
$y=(\frac{m-1}2,\frac{m-1}2)$. Hence, since $diam(C_{2k+1})=k<m-1$,
we have that  $\Lambda\times\{\frac{m-1}2\}\subseteq I^2[u,v]$. This fact together with  Lemma~\ref{lem:rodaja}, allows us to derive  that 
$\{ u,v\}$ is  a hull set of $C_{2k+1}\boxtimes P_m$.

\begin{figure}[h]
\begin{center}
  \includegraphics[width=0.8\textwidth]{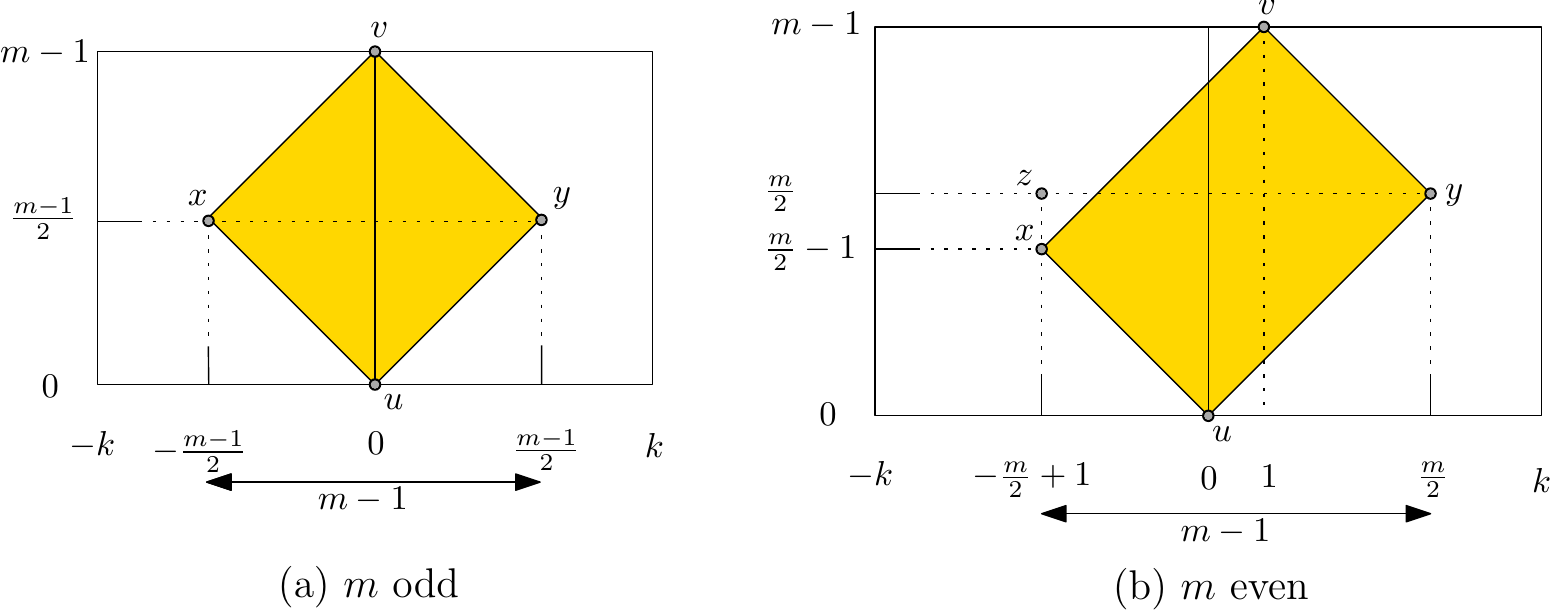}
 \caption{In both cases, $k\le m-2$.}
  \label{fig.caminoHull2}
\end{center}
\end{figure}

Next, assume that $k\le m-2$, $m$ is even and  take the pair of vertices
$u=(0,0)$ and $v=(1,m-1)\}$.
Notice that, as shown in Figure \ref{fig.caminoHull2} (b), $I[u,v]$ is 
the set of points belonging to  the rectangle determined by
$\{u,v,x,y\}$, where $x=(-\frac{m}2+1,\frac{m}2-1)$ and $y=(\frac{m}2,\frac{m}2)$.
Hence, since $diam(C_{2k+1})=k<m-1=\frac m2+(\frac m2 -1)$, we have that 
$\Lambda\times\{\frac{m}2\}\setminus\{z=(-\frac{m}2+1,\frac{m}2)\}\subseteq 
I^2[u,v]$,
i.e., $\Lambda\times\{\frac{m}2\}\subseteq I^3[u,v]$.
This fact together with  Lemma~\ref{lem:rodaja}, allows us to derive  
that  $\{ u,v\}$ is  a hull set of $C_{2k+1}\boxtimes P_m$.

Finally, assume that $k\ge m-1$ and take an arbitrary 2-vertex set $\{u,v\}\subset V(C_{2k+1}\boxtimes C_{2k+1})$. We may assume wlog that $u=(0,h)$ and $v=(a,h')$, where $0\le a\le k$ and $0\le h\le h'$. We distinguish two cases.

{\bf Case 1.} $d(u,v)=\max \{ a, h'-h \}=a>0$: The path $\rho$ of $C_{2k+1}$ whose vertex set is $V(\rho)=\{0,1,2,\dots ,a\}$ is the unique $0-a$ geodesic. Hence, according to Lemma \ref{lem.projmax}, 
$I[u,v]$ is the subset of points of $\mathcal{P}$ lying in the rectangle $\mathcal{R}$ determined by the four lines passing through either $u$ or $v$, of slopes $\pm 1$. Note that, as shown in  Figure \ref{fig.cicloCamino} (a), this rectangle is inside the square of side length $a$ determined by the four vertices of $\mathcal{R}$. 
This fact, together with Lemma \ref{lem.projmax}, implies that for any pair of vertices $u',v'\in I[u,v]=\mathcal{R}$, the set $I[u',v']$ is the rectangle $\mathcal{R'}$ contained in $\mathcal{R}$, determined by the four lines passing through either $u'$ or $v'$, of slopes $\pm 1$. This means that $I[u,v]$ is a proper convex subset of $C_{2k+1}\boxtimes P_{m}$, and thus $\{u,v\}$ is not geodetic.

\begin{figure}[h]
\begin{center}
  \includegraphics[width=0.9\textwidth]{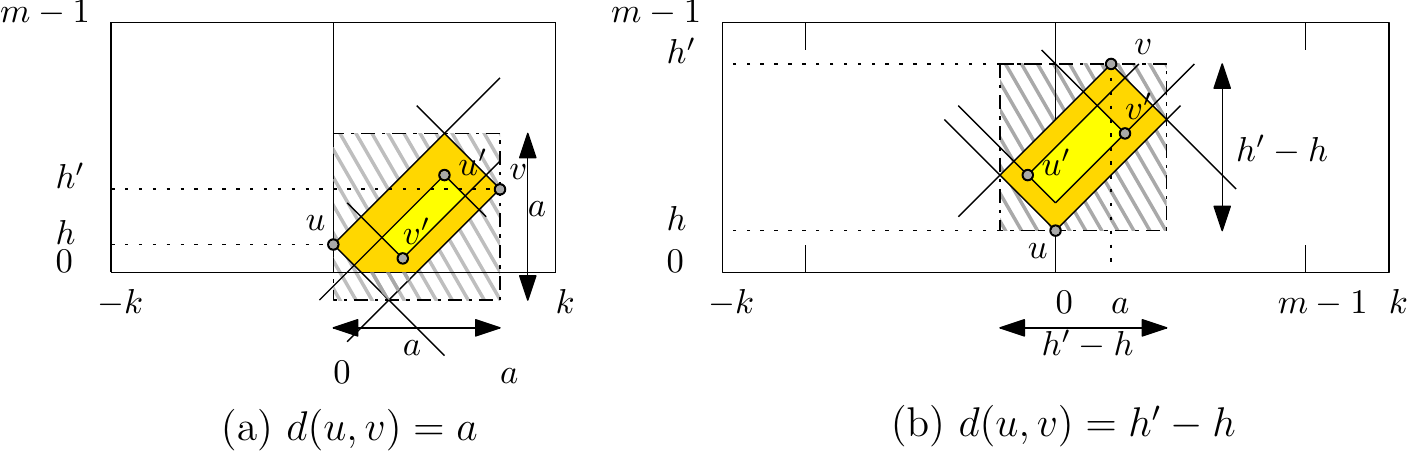}
 \caption{In both cases, $\{u=(0,h),v=(a,h')\}\subset V(C_{2k+1}\boxtimes P_{m})$.}
  \label{fig.cicloCamino}
\end{center}
\end{figure}

{\bf Case 2.} $d(u,v)=\max \{ a, h'-h \}=h'-h>0$: As shown in  Figure \ref{fig.cicloCamino} (b), $I[u,v]$
consists of all  points of $\mathcal{P}$ lying in the rectangle $\mathcal{R}$ determined by the four lines passing through either $u$ or $v$, of slopes $\pm 1$, And   it is inside the square of side length $h'-h\le m-1\le k$ determined by the four vertices of $\mathcal{R}$. Hence,
reasoning as in the preceding case and having in mind that $h'-h\le k$,
we derive that  $I[u,v]$ is a proper convex subset of $C_{2k+1}\boxtimes P_{m}$, and thus $\{u,v\}$ is not geodetic.
\end{proof}

The last strong product graphs we have studied is the so-called family of strong toruses, i.e., the strong product of two  cycles. As an immediate  consequence of Proposition~\ref{prop.gno3} and Theorem~\ref{teo.geodetico}, the following results are derived.

\begin{prop}\label{toruses} Let $m,n$ be two  integers such that $4\le\min\{m,n\}$. 
\begin{itemize} 
\item[(i)] if both $m$ and $n$ are even, then $g(C_m\boxtimes C_n)=4$
\item[(ii)]  if $mn$ is even, then $4\le g(C_m\boxtimes C_n)\le 6$.
\end{itemize}
\end{prop}

In addition, we have been able to obtain a number of further results for the geodetic number, involving odd cycles.

\begin{figure}[htb]
\begin{center}
\includegraphics[width=0.7\textwidth]{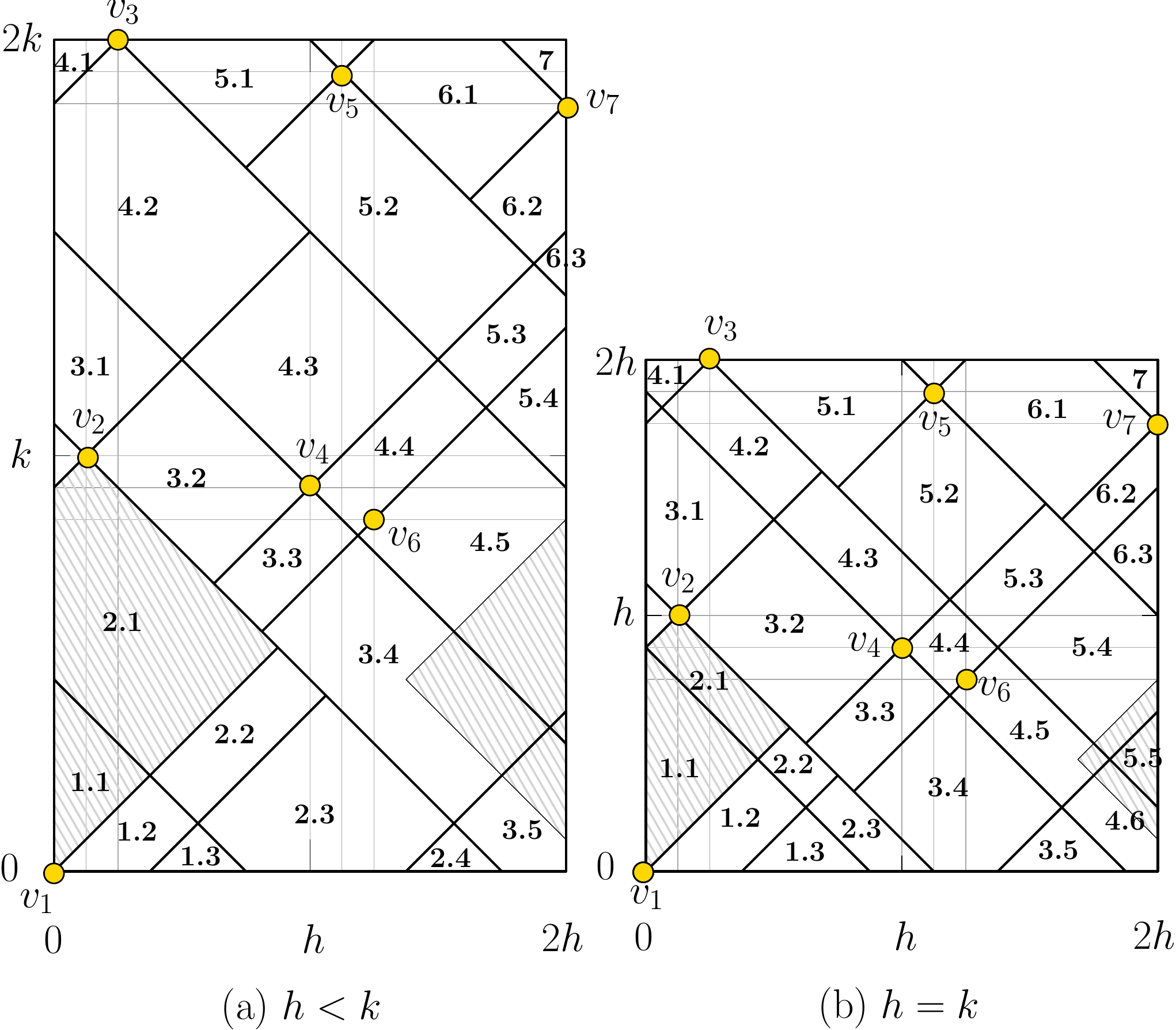}
 \caption{Partition of   $V(C_{2h+1}\boxtimes C_{2k+1})$. The dashed  region is $I[v_1,v_2]$.}
  \label{geoCiclosm7}
\end{center}
\end{figure}

\begin{prop} Let $h,k,n$ be integers such that $2\le\min\{h,k\}$ and $4\le n$. 
\begin{itemize} 

\item[(i)] If $2h \le k$, then $g(C_{2h+1}\boxtimes C_{2k})= 4$.

\item[(ii)] If $5(2h-1) \le n$, then $g(C_{2h+1}\boxtimes C_{n})\le 5$.

\item[(iii)] If $2h+1 \le k$, then $g(C_{2h+1}\boxtimes C_{2k+1})\le 6$.

\item[(iv)] If $3\le h \le k$, then $5\le g(C_{2h+1}\boxtimes C_{2k+1})\le7$.

\end{itemize}

\end{prop}
\begin{proof} (i) It is  straightforward  to verify that $V(C_{2h+1}\boxtimes C_{2k})=I[S_1]=I[u_1,u_3]\cup I[u_2,u_4]$, where         $S_1=\{u_1,u_2,u_3,u_4\}$, $u_1=(0,0)$, $u_2=(h,h)$, $u_3=(0,k)$ and $u_4=(h,h+k-1)$.

(ii) It is  straightforward  to verify that $V(C_{2h+1}\boxtimes C_{n})=I[S_2]=I[w_1,w_3]\cup I[w_1,w_4]\cup I[w_2,w_4]\cup I[w_2,w_5]\cup I[w_3,w_5]$, where 
    $S_2=\{w_1,w_2,w_3,w_4,w_5\}$, $t=\lfloor\frac{n}{5}\rfloor$,  $w_1=(h,0)$, $w_2=(h,t)$, $w_3=(h,2t)$, $w_4=(h,3t)$ and $w_5=(h,4t)$. 

(iii) It is  straightforward  to verify that $V(C_{2h+1}\boxtimes C_{2k+1})=I[S_3]=I[u_1,u_3]\cup I[u_2,u_4]\cup I[u_2,u_6]\cup I[u_3,u_5]$, where 
$S_3=\{u_1,u_2,u_3,u_4,u_5,u_6\}$, $u_5=(0,2k)$ and $u_6=(h,h+k)$.

(iv) The lower bound is a direct consequence of Lemma \ref{vertincycles} and Proposition \ref{ultima}. The upper bound is obtained as a consequence of 
the following claim:
$V(C_{2h+1}\boxtimes C_{2k+1})=I[S_4]=I[v_1,v_2]\cup I[v_1,v_5]\cup I[v_1,v_7]\cup I[v_2,v_3]\cup I[v_2,v_4]\cup I[v_2,v_6]\cup I[v_3,v_4]\cup I[v_3,v_5]\cup I[v_4,v_5]\cup I[v_5,v_6]\cup I[v_5,v_7]\cup I[v_6,v_7]$, where $S_4=\{v_1,v_2,v_3,v_4,v_5,v_6,v_7\}$, $v_1=(0,0)$, $v_2=(1,k)$, $v_3=(2,2k)$, $v_4=(h,k-1)$, $v_5=(h+1,2k-1)$, $v_6=(h+2,k-2)$ and $v_7=(2h,2k-2)$.

To prove this claim, let us first partition the vertex set of $C_{2h+1}\boxtimes C_{2k+1}$ as shown in Figure \ref{geoCiclosm7}, according to the following facts:

{\scriptsize
\begin{itemize}

\item[(1)] $0 \le i+j \le h-1$ and $\left\{\begin{array}{lcll}
                                                (1.1)  & i-j \le 0 &  \Rightarrow &  (i,j)\in I[v_1,v_2] \\
                                                (1.2)  & 0 < i-j \le 3 &  \Rightarrow &  (i,j)\in I[v_3,v_4] \\
                                                (1.3)  & 3 < i-j  &  \Rightarrow &  (i,j)\in I[v_3,v_5]
                                              \end{array}
                                            \right.$
\item[(2)] $h-1 < i+j < k+1$ and $\left\{\begin{array}{lcll}
                                                (2.0)  & i=0, j=k &  \Rightarrow &  (i,j)\in I[v_2,v_6] \\
                                                (2.1)  & 1-k \le i-j \le 0 &  \Rightarrow &  (i,j)\in I[v_1,v_2] \\
                                                (2.2)  & 0 < i-j \le 3 &  \Rightarrow &  (i,j)\in I[v_3,v_4] \\
                                                (2.3)  & 3 < i-j \le h+3 &  \Rightarrow &  (i,j)\in I[v_3,v_5] \\
                                                (2.4)  & h+3 < i-j  &  \Rightarrow &  (i,j)\in I[v_1,v_5] 
                                              \end{array}
                                            \right.$

\item[(3)] $k+1 \le i+j \le h+k-1$ and $\left\{\begin{array}{lcll}
                                              (3.1)  & i-j \le 1-k &  \Rightarrow &  (i,j)\in I[v_2,v_3] \\
                                              (3.2)  & 1-k < i-j \le 1+h-k &  \Rightarrow &  (i,j)\in I[v_2,v_4] \\
                                              (3.3)  & 1+h-k < i-j < 4+h-k &  \Rightarrow &  (i,j)\in I[v_3,v_4] \\
                                              (3.4)  & 4+h-k \le i-j \le h+3 &  \Rightarrow &  (i,j)\in I[v_5,v_6] \\
                                              (3.5)  & h+3 < i-j & \Rightarrow &  (i,j)\in I[v_1,v_2]\cup I[v_1,v_5]
                                              \end{array}
                                            \right.$

\item[(4)] $h+k-1 < i+j < 2k+2$ and $\left\{\begin{array}{lcll}
                                            (4.1)  & i-j < 2-2k &  \Rightarrow &  (i,j)\in I[v_1,v_7] \\
                                            (4.2)  & 2-2k \le i-j < 1-k &  \Rightarrow &  (i,j)\in I[v_2,v_3] \\
                                            (4.3)  & 1-k \le i-j \le 1+h-k &  \Rightarrow &  (i,j)\in I[v_4,v_5] \\
                                            (4.4)  & 1+h-k < i-j \le 4+h-k &  \Rightarrow &  (i,j)\in I[v_6,v_7] \\
                                            (4.5)  & 4+h-k <i-j \le h+3 &  \Rightarrow &  (i,j)\in I[v_2,v_6] \\
                                            (4.6)  & h+3 < i-j &  \Rightarrow &  (i,j)\in I[v_1,v_5]
                                              \end{array}
                                            \right.$

\item[(5)] $2k+2 \le i+j \le h+2k$ and $\left\{\begin{array}{lcll}
                                              (5.1)  & i-j \le 2+h-2k &  \Rightarrow &  (i,j)\in I[v_3,v_5] \\
                                              (5.2)  & 2+h-2k < i-j \le 1+h-k &  \Rightarrow &  (i,j)\in I[v_4,v_5] \\
                                              (5.3)  & 1+h-k < i-j \le 4+h-k &  \Rightarrow &  (i,j)\in I[v_6,v_7] \\
                                              (5.4)  & 4+h-k < i-j \le h+3 &  \Rightarrow &  (i,j)\in I[v_2,v_6] \\
                                              (5.5)  & h+3 < i-j & \Rightarrow &  (i,j)\in I[v_1,v_2]
                                              \end{array}
                                            \right.$

\item[(6)] $h+2k \le i+j \le 2h+2k-2$ and $\left\{\begin{array}{lcll}
                                            (6.0)  & i-=h+1, j=2k &  \Rightarrow &  (i,j)\in I[v_5,v_6] \\
                                            (6.1)  & 2+h-2k \le i-j \le 2+2h-2k &  \Rightarrow &  (i,j)\in I[v_5,v_7] \\
                                            (6.2)  & 2+2h-2k < i-j \le 4+h-k &  \Rightarrow &  (i,j)\in I[v_6,v_7] \\
                                            (6.3)  & 4+h-k < i-j &  \Rightarrow &  (i,j)\in I[v_2,v_6]
                                            \end{array}
                                            \right.$

\item[(7)] $2h+2k-2 < i+j \le 2h+2k $ $\Rightarrow$ $(i,j)\in I[v_1,v_7]$

\end{itemize}
}

It is straightforward to see that each of the regions described above is contained in at least one of the 12 mentioned  closed intervals. As a matter of example, notice that if $h=k$, then regions (1.1), (2.1) and (5.5) are completely covered by $I[v_1,v_2]$, since:

\vspace{.3cm}$(i,j)\in I[v_1,v_2] \Leftrightarrow$
$\left\{\begin{array}{ll}
0\le i \le 2h \\
0\le j \le 2k
\end{array}\right.$   and either 
$\left\{\begin{array}{ll}
1-k\le i-j \le 0 \\
0\le i+j \le k+1
\end{array}\right.$ or
$\left\{\begin{array}{ll}
2h-k+2\le i-j  \\
2h+1 \le i+j 
\end{array}\right.$
\end{proof}

Let us remark that all  bounds  presented in the last two propositions can not be improved, as it is shown in Table \ref{t7}, which contains  the geodetic number of some strong product graphs of the form $C_5\boxtimes C_n$ computationally obtained.

\begin{table}[htbp]
\begin{center}\small
\begin{tabular}{|c||c|c|c|c|c|c|}
  \hline
  
$C_5\boxtimes C_n$  & $C_5\boxtimes C_4$ & $C_5\boxtimes C_5$  & $C_5\boxtimes C_6$ & $C_5\boxtimes C_7$ & $C_5\boxtimes C_8$ & $C_5\boxtimes C_9$\\
\hline
\hline
  
$g(C_5\boxtimes C_n)$  & $5$ & $5$  & $6$ & $7$ & $4$ & $6$\\
\hline
\hline

\end{tabular}
\end{center}
\vspace{-.4cm}\caption{ Geodetic number of some strong product graphs of the form $C_5\boxtimes C_n$.\label{t7}}
\end{table}

\begin{figure}[htb]
\begin{center}
\includegraphics[width=0.6\textwidth]{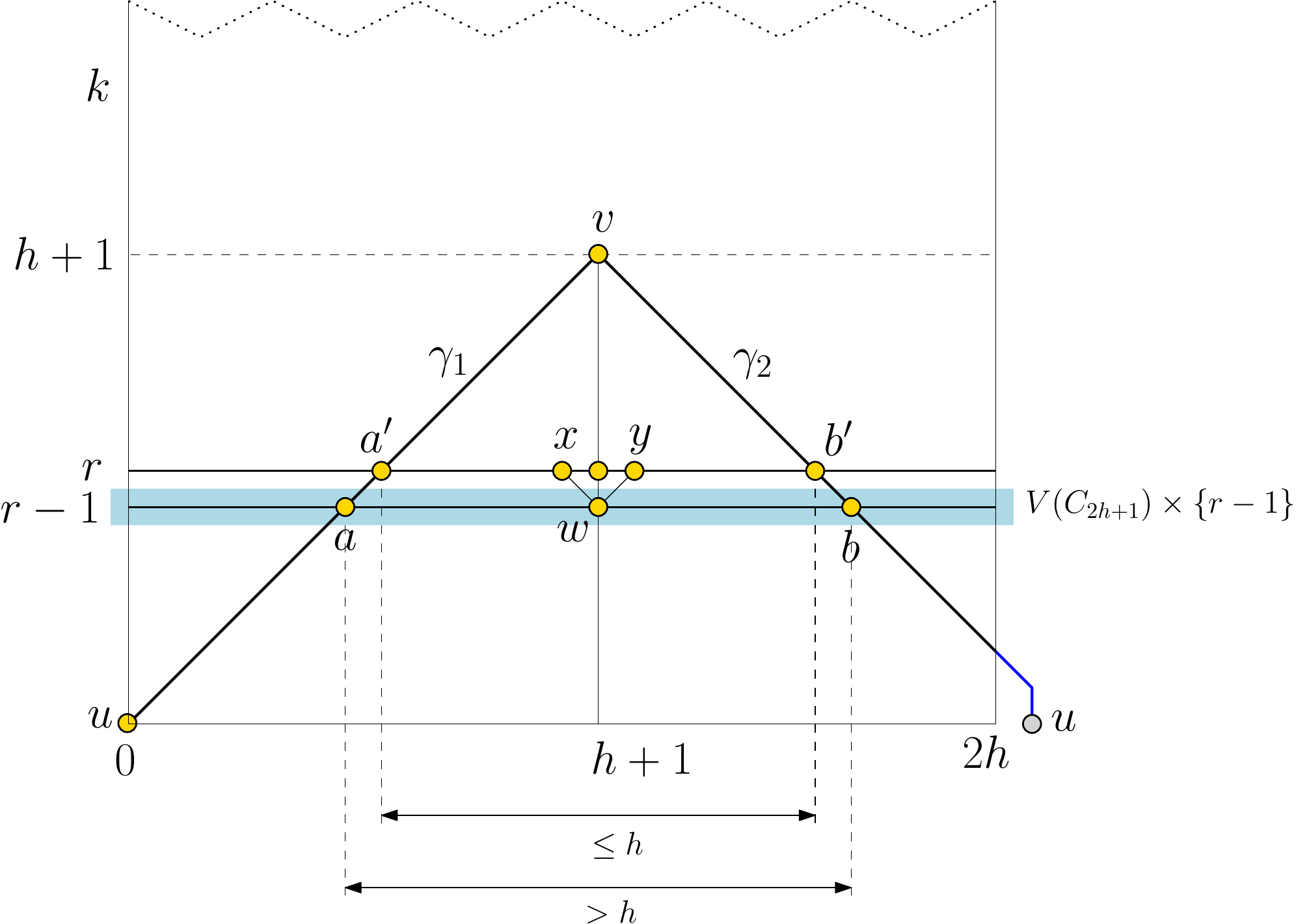}
 \caption{  $\gamma_1$ and $\gamma_2$ are two $u-v$ geodesics in $C_{2h+1}\boxtimes C_{2k+1}$, whenever  $h<k$.}
  \label{fig.2h3}
\end{center}
\end{figure}

\begin{figure}[htb]
  \begin{center}
  \includegraphics[width=0.37\textwidth]{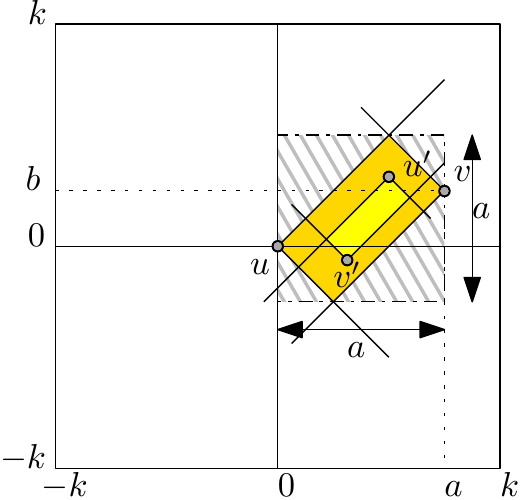}
  \caption{$u=(0,0)$ and $v=(a,b)$ are two vertices of the strong torus $C_{2k+1}\boxtimes C_{2k+1}$}
  \label{fig.cicimp}
  \end{center}
\end{figure}

\begin{prop} If $h,k$ are  integers such that $2\le h\le k$, then $h(C_{2h+1}\boxtimes C_{2k+1})=\left\{
                                              \begin{array}{ll}
                                                2, & \hbox{if }h< k; \\
                                                3, & \hbox{if }h=k.
                                              \end{array}
                                            \right.$
\end{prop}
\begin{proof} 
Certainly, $2\le h(C_{2k+1}\boxtimes C_{2k+1})\le 3$, being the upper bound  a corollary of Theorem~\ref{teo.hull2}, whereas the lower bound is derived from the fact that  $h(G)=1$ if and only if $G=K_1$.

Suppose next that $h<k$ and consider the set $S=\{u,v\}$, where $u=(0,0)$ and $v=(h+1,h+1)$ (see Figure \ref{fig.2h3}).  Since $d(u,v)=h+1$, we have  that both
$\gamma_1: (0,0)(1,1)\ldots(h,h)(h+1,h+1)$ and $\gamma_2:(0,0)(0,1)(2h,2)\ldots(h+2,h)(h+1,h+1)$ are  $u-v$ geodesic. Observe that given $c=(i,i)\in V(\gamma_1)$ and $d=(2h+2-i,i)\in V(\gamma_2)$, $d(c,d)=2h+2-2i$ if and only if 
$2h+2-2i\le h$, i.e., if and only if $\lceil\frac{h+2}{2}\rceil\le i$. At this point, we claim that if $r=\lceil\frac{h+2}{2}\rceil$, then  $V(C_{2h+1})\times\{r-1\}\subseteq I^4[S]$, which according to Lemma~\ref{lem:rodaja}, is enough to end the proof of the case $h<k$.

To show that this claim is true, consider the vertices $a=(r-1,r-1)$, $w=(h+1,r-1)$, $b=(2h-r+3,r-1)$, $a'=(r,r)$, $x=(h,r)$, $y=(h+2,r)$, $b'=(2h-r+2,r)$ and observe:


$\ast$ $\{a,b,a',b'\}\subset I[S]$, since $\{a,a'\}\subset V(\gamma_1)$ and $\{b,b'\}\subset V(\gamma_2)$.

$\ast$ $\{x,y\}\subset I[a',b']\subseteq I^2[S]$, since $r<h<h+2<2h-r+2$ and $d(a',b')=2h-2r+2\le h$.

$\ast$ $w\in I[x,y] \subseteq I^3[S]$, since $d(x,w)=d(w,y)=1$ and $d(x,y)=2$.


Finally, we show that for every  $i\in V(C_{2h+1})$, the vertex $z=(i,r-1)\in I^4[S]$:

$\bullet$ If $0\le i \le  r-1$, then $z\in I[a,b]\subset I^2[S]$, since $(2h-r+3)-(r-1)=2h-2r+4\ge h+1$.

$\bullet$ If $r-1 \le i \le  h+1$, then $z\in I[a,w]\subset I^4[S]$, since  $(h+1)-(r-1)=h-r+2\le h$.

$\bullet$ If $h+1 \le i \le 2h-r+3$, then $z\in I[w,b]\subset I^4[S]$, since $(2h-r+3)-(h+1)=h-r+2\le h$.

$\bullet$ If $2h-r+3 \le i \le 2h$, then $z\in I[a,b]\subset I^2[S]$, since $(2h-r+3)-(r-1)\ge h+1$.

Now, suppose that $h=k$ and take an arbitrary 2-vertex set $\{u,v\}\subset V(C_{2k+1}\boxtimes C_{2k+1})$. By symmetry reasons, we label the vertex set of $C_{2k+1}$ as follows: $V(C_{2k+1})=\Lambda=\{ -k, \dots ,-1,0,1,\dots ,k\}$.

We may assume wlog that $u=(0,0)$ and $v=(a,b)$, where $0\le b\le a\le k$. Observe that $d(u,v)=\max \{ a,b \}=a\le k$, and  that the path $\rho$ of $C_{2k+1}$ whose vertex set is $V(\rho)=\{0,1,2,\dots ,a\}$ is the unique $0-a$ geodesic. Hence, according to Lemma \ref{lem.projmax}, 
$$I[u,v]=\{ (i,j) : 0\le i\le a \, , \, |j|\le i\, , \, |j-b|\le a-i\}=\{ (i,j) :  0 \le i-j \le a-b  \, , \, 0 \le i+j \le a+b\}$$ 
 
In other words, if we identify $V(C_{2k+1}\boxtimes C_{2k+1})$ with the grid $\mathcal{P}=\Lambda\times\Lambda$ of points of the discrete plane $\mathbb{Z}^2$, then $I[u,v]$ is the subset of points of $\mathcal{P}$ lying in the rectangle $\mathcal{R}$ determined by the four lines passing through either $u$ or $v$, of slopes $\pm 1$. Note that, as shown in  Figure \ref{fig.cicimp}, this rectangle is inside the square of side length $a$ determined by the four vertices of $\mathcal{R}$. 
This fact, together with Lemma \ref{lem.projmax}, implies that for any pair of vertices $u',v'\in I[u,v]=\mathcal{R}$, the set $I[u',v']$ is the rectangle $\mathcal{R'}$ contained in $\mathcal{R}$, determined by the four lines passing through either $u'$ or $v'$, of slopes $\pm 1$. This means that $I[u,v]$ is a proper convex subset of $C_{2k+1}\boxtimes C_{2k+1}$, and thus $\{u,v\}$ is not geodetic.
\end{proof}

We conclude this section by showing in Table \ref{t5} all the results obtained for strong product graphs of the form $G\boxtimes C_n$.

\begin{table}[htbp]
\begin{center}\small
\begin{tabular}{|c||c|c|}
  \hline
  
  $G\boxtimes C_n$  & {\bf $g(G\boxtimes C_n)$} & {\bf $h(G\boxtimes C_n)$} \\
\hline
\hline
  
$P_m\boxtimes C_n$  

&  $\left\{
  \begin{array}{ll}
    4, & \hbox{ if \emph{n} is even;} \\
    { 5,6,} & { \hbox{ if \emph{n} is odd.}}
  \end{array}
\right.$   

&  $\left\{
  \begin{array}{ll}
     {3,} &  {\hbox{if $n=2r+1$ odd and $m<r+2$};} \\
    2, & \hbox{ otherwise.}
  \end{array}
\right.$  \\  \hline

$K_m\boxtimes C_n$  

&   $\left\{
  \begin{array}{ll}
    4, & \hbox{ if \emph{n} is even;} \\
    5, & \hbox{ if \emph{n} is odd.}
  \end{array}
\right.$   

&   $\left\{
  \begin{array}{ll}
    2, & \hbox{ if \emph{n} is even;} \\
     3, &   \hbox{ if \emph{n} is odd.}
  \end{array}
\right.$  \\  \hline

$C_m\boxtimes C_n$  

&   $\left\{
  \begin{array}{ll}
    4, & \hbox{ if \emph{m} and \emph{n} are even;} \\
    4-6, & \hbox{ if \emph{m} is even and \emph{n} is odd;} \\
    {5-7,} & { \hbox{ if \emph{m} and \emph{n} are odd.}}
  \end{array}
\right.$   

&   $\left\{
   \begin{array}{ll}
    3, & \hbox{if $m=n$ is odd;} \\
    2, & \hbox{  otherwise.}
  \end{array}
\right.$ \\  \hline
\hline
\end{tabular}
\end{center}
\vspace{-.4cm}\caption{ Geodetic and hull numbers of some strong product graphs of the form $G\boxtimes C_n$.\label{t5}}
\end{table}

%
%
%
%
%
%

\bibliographystyle{elsarticle-num}



\end{document}